# On the Homotopy Type of the Fredholm Lagrangian-Grassmannian

J. C. C. Eidam, P. Piccione*

February 24, 2005

**Abstract**

In this article, using the results of [1], we show that the Fredholm Lagrangian-Grassmannian $\mathscr{F}_{\mathsf{L}_0}(\Lambda)$ has the homotopy type of the quotient $\mathrm{U}(\infty)/\mathrm{O}(\infty)$. As a corollary, we obtain that the Maslov index with respect to a lagrangian is a isomorphism between the fundamental group of $\mathscr{F}_{\mathsf{L}_0}(\Lambda)$ and the integers. We remark that the fundamental group $\pi_1(\mathscr{F}_{\mathsf{L}_0}(\Lambda))$ was computed by B. Booss-Bavnbek and K. Furutani in [2], [4], with considerable effort, but the method works only for the $\pi_1$ case.

# 1 Introduction

Let $\mathsf{H}$ be a separable infinite dimensional Hilbert space endowed with a sympletic structure, that is, there exists a bounded operator $J$ in $\mathsf{H}$ such that $J^2 = -I$ and $J^* = -J$ and a sympletic form $\omega$ given by $\omega(\cdot,\cdot) = \langle J\cdot,\cdot\rangle$. The space $\Lambda$ of all lagrangian subspaces in $\mathsf{H}$ admits a structure of real-analytic contractible Banach manifold, so uninteresting from the topological viewpoint. To obtain a space with nontrivial topology, consider $\mathsf{L}_0 \in \Lambda$ and the space $\mathscr{F}_{\mathsf{L}_0}(\Lambda)$ of all $\mathsf{L} \in \Lambda$ which forms a Fredholm pair with $\mathsf{L}_0$. It is a open set in $\Lambda$, so it is also a analytic Banach manifold.

In this article, we prove that $\mathscr{F}_{\mathsf{L}_0}(\Lambda)$ is homotopy equivalent with the space of compact perturbations of a fixed lagrangian; this last space is homotopy equivalent with $\mathrm{U}(\infty)/\mathrm{O}(\infty)$.

# 2 Proof of the result

Let us fix some notation. Denote by $P_{\mathsf{M}}$ the orthogonal projection over a closed subspace $\mathsf{M} \subset \mathsf{H}$. Denote by $\mathrm{GL}(\mathsf{H})$ the group of all invertible operators in $\mathsf{H}$ and by $\mathrm{O}(\mathsf{H})$ the group of unitary operators in $\mathsf{H}$. If $\mathsf{H}$ is a sympletic space with complex structure $J$, the group $\mathrm{U}(\mathsf{H};J)$ is formed by elements of $\mathrm{O}(\mathsf{H})$ which commute with $J$. We define $\mathrm{GL}_c(\mathsf{H}) = \mathrm{GL}(\mathsf{H}) \cap (I + \mathscr{L}_c(\mathsf{H}))$, $\mathrm{O}_c(\mathsf{H}) = \mathrm{O}(\mathsf{H}) \cap (I + \mathscr{L}_c(\mathsf{H}))$ and $\mathrm{U}_c(\mathsf{H};J) =$

*Research supported by Fapesp, São Paulo, Brasil, grants 01/00046-3 and 02/02528-8. The second author is partially sponsored by Cnpq, Brasil.

$\mathrm{U}(\mathsf{H};J)\cap(I+\mathscr{L}_c(\mathsf{H}))$, where $\mathscr{L}_c(\mathsf{H})$ denotes the space of compact operators in $\mathsf{H}$. Elements of $I+\mathscr{L}_c(\mathsf{H})$ are called *compact perturbations of identity*.

We begin by studying some properties of compact perturbations of lagrangians. Given $\mathsf{L},\mathsf{L}'\in\Lambda$, $\mathsf{L}'$ is a compact perturbation of $\mathsf{L}$ if $P_{\mathsf{L}'}-P_{\mathsf{L}}$ is compact. This fact is denoted by $\mathsf{L}'\sim\mathsf{L}$. Obviously, $\sim$ is a equivalence relation in $\Lambda$. The quotient space of this relation is called *essential grassmannian-lagrangian* and is denoted by $\Lambda_e$. It is endowed with the quotient topology induced by the canonical map $\pi:\Lambda\to\Lambda_e$. For $\mathsf{L}_0\in\Lambda$, we define $\mathcal{G}_c(\mathsf{L}_0)$ as the equivalence class of $\mathsf{L}_0$. The next results determine the topology of $\mathcal{G}_c(\mathsf{L}_0)$.

Lemma 2.1 *The group action*

$$\begin{array}{rcl}\mathrm{U}_c(\mathsf{H};J)\times\mathcal{G}_c(\mathsf{L}_0) & \to & \mathcal{G}_c(\mathsf{L}_0)\\ (U,\mathsf{L}) & \mapsto & U(\mathsf{L})\end{array}\tag{2.1}$$

*is transitive.*

Proof. Let $\mathsf{L}\in\mathcal{G}_c(\mathsf{L}_0)$; we must find $U\in\mathrm{U}_c(\mathsf{H};J)$ such that $U(\mathsf{L}_0)=\mathsf{L}$. Define $T_1=I+(P_{\mathsf{L}}-P_{\mathsf{L}_0})|_{\mathsf{L}_0}=P_{\mathsf{L}}|_{\mathsf{L}_0}:\mathsf{L}_0\to\mathsf{L}$; since $T_1$ is a Fredholm operator with $\operatorname{ind}T_1=0$, there exists $K:\mathsf{L}_0\to\mathsf{L}$ of finite rank such that $T_2=T_1+K:\mathsf{L}_0\to\mathsf{L}$ is a linear isomorphism. Putting

$$T=T_2P_{\mathsf{L}_0}-JT_2JP_{\mathsf{L}_0^\perp}=I+((P_{\mathsf{L}}-P_{\mathsf{L}_0})+K)P_{\mathsf{L}_0}-J((P_{\mathsf{L}}-P_{\mathsf{L}_0})+K)JP_{\mathsf{L}_0^\perp}\,,$$

we have that $T\in\mathrm{GL}_c(\mathsf{H})$, $TJ=JT$ and $T(\mathsf{L}_0)=\mathsf{L}$. This last equality implies that $T^*(\mathsf{L}^\perp)=\mathsf{L}_0^\perp$, so, since $T^*J=JT^*$, we have

$$T^*(\mathsf{L})=T^*J(\mathsf{L}^\perp)=JT^*(\mathsf{L}^\perp)=J(\mathsf{L}_0^\perp)=\mathsf{L}_0\,.$$

In particular, $(T^*T)(\mathsf{L}_0)=\mathsf{L}_0$. Let $S$ be the unique positive square root for $T^*T$. It is easy (see for example [10, Chapter VII]) that $S$ can be given as the limit of the sequence $\{S_n\}$ given by $S_0=0$ and

$$S_{n+1}=\frac{1}{2}(T^*T+S_n^2)\,,\tag{2.2}$$

for $n\geq 0$. Since $T$ is a compact perturbation of the identity, we see, by formula (2.2) and induction that $S$ is a compact perturbation of the identity, $SJ=JS$ and $S(\mathsf{L}_0)=\mathsf{L}_0$. Taking $T=US$ the polar decomposition of $T$, with $U=TS^{-1}$, we see that $U(\mathsf{L}_0)=\mathsf{L}$ and $U\in\mathrm{U}_c(\mathsf{H};J)$, as desired. ∎

Proposition 2.2 *The space $\mathcal{G}_c(\mathsf{L}_0)$ is homotopy equivalent to* $\mathrm{U}(\infty)/\mathrm{O}(\infty)$.

Proof. Consider the differentiable action (2.1). Since the isotropy group at $\mathsf{L}_0$ is $\{U\in\mathrm{U}_c(\mathsf{H};J):U(\mathsf{L}_0)=\mathsf{L}_0\}\approx\mathrm{O}_c(\mathsf{L}_0)$, we conclude that there is a diffeomorfism $\mathcal{G}_c(\mathsf{L}_0)\approx\mathrm{U}_c(\mathsf{H};J)/\mathrm{O}_c(\mathsf{L}_0)$. Since the spaces $\mathrm{U}_c(\mathsf{H};J)$ and $\mathrm{O}_c(\mathsf{L}_0)$ are homotopy equivalent with $\mathrm{U}(\infty)$ and $\mathrm{O}(\infty)$, respectively, (see [7]), the result follows. ∎

If $U\in\mathrm{U}(\mathsf{H};J)$ and $\mathsf{L}\in\Lambda$ then $P_{U(\mathsf{L})}=UP_{\mathsf{L}}U^*$. This fact is important to prove the next theorem, which describes the topology of $\Lambda_e$.

Theorem 2.3 *The canonical map $\pi:\Lambda\to\Lambda_e$ defines a fibre bundle with typical fiber $\mathcal{G}_c(\mathsf{L}_0)$.*

Proof. Let $\mathsf{L}_0 \in \Lambda$ and consider the Calkin algebra $\mathcal{C}(\mathsf{H}) = \mathscr{L}(\mathsf{H})/\mathscr{L}_c(\mathsf{H})$, where $\mathscr{L}(\mathsf{H})$ ($\mathscr{L}_c(\mathsf{H})$) denotes the space of bounded (compact) operators in $\mathsf{H}$. As proved in [1, Proposition 1.2], if $\mathsf{U}(\mathcal{C}(\mathsf{H}))$ denotes the unitary elements of $\mathcal{C}(\mathsf{H})$ (that is, $u \in \mathsf{U}(\mathcal{C}(\mathsf{H}))$ if and only if $u^*u = uu^* = 1$), then the map

$$\begin{array}{rcl} \mathsf{U}(\mathcal{C}(\mathsf{H})) & \to & \mathsf{U}(\mathcal{C}(\mathsf{H})) \cdot [P_{\mathsf{L}_0}] \\ u & \mapsto & u[P_{\mathsf{L}_0}]u^* \end{array}$$

is a analytic principal bundle with structural group $\{u \in \mathsf{U}(\mathcal{C}(\mathsf{H})) \,:\, u[P_{\mathsf{L}_0}]u^* = [P_{\mathsf{L}_0}]\}$, where $[P_{\mathsf{L}_0}]$ denotes the image of $P_{\mathsf{L}_0}$ by the canonical map $p : \mathscr{L}(\mathsf{H}) \to \mathcal{C}(\mathsf{H})$. By restricting the structural group to the closed subgroup

$$G = p(\{U \in \mathsf{U}(\mathsf{H}; J) \,:\, U(\mathsf{L}_0) = \mathsf{L}_0\}),$$

we obtain a principal bundle

$$\begin{array}{rcl} \Phi : p(\mathsf{U}(\mathsf{H}; J)) & \to & p(\mathsf{U}(\mathsf{H}; J)) \cdot [P_{\mathsf{L}_0}] \\ u & \mapsto & u[P_{\mathsf{L}_0}]u^* \end{array} \tag{2.3}$$

with $G$ as structural group. The map (2.3) can be rewritten as

$$\begin{array}{rcl} \Phi : p(\mathsf{U}(\mathsf{H}; J)) & \to & p(\mathsf{U}(\mathsf{H}; J)) \cdot [P_{\mathsf{L}_0}] \approx \Lambda_e \\ u & \mapsto & [P_{U(\mathsf{L}_0)}], \end{array}$$

where $U \in \mathsf{U}(\mathsf{H}; J)$ is such that $p(U) = u$. Since we can obtain local sections for the map $\Phi$ and the orbit $p(\mathsf{U}(\mathsf{H}; J)) \cdot [P_{\mathsf{L}_0}]$ can be identified with the space $\Lambda_e$, there exists a neighbourhood $\mathcal{U} \subset \Lambda_e$ of $\pi(\mathsf{L}_0)$ and a map $s' : \mathcal{U} \to p(\mathsf{U}(\mathsf{H}; J))$ such that $\Phi \circ s' = I$. Now, again by the results of [1, Proposition 1.7], the map $p|_{\mathsf{U}(\mathsf{H};J)} : \mathsf{U}(\mathsf{H}; J) \to p(\mathsf{U}(\mathsf{H}; J))$ admit local sections around all points of $p(\mathsf{U}(\mathsf{H}; J))$, so, we obtain a neighbourhood $\mathcal{V} \subset p(\mathsf{U}(\mathsf{H}; J))$ of $I$ and a map $s'' : \mathcal{V} \to \mathsf{U}(\mathsf{H}; J)$ such that $p \circ s'' = I$. Putting $s = s'' \circ s' : \mathcal{U} \to \mathsf{U}(\mathsf{H}; J)$, the map

$$\begin{array}{rcl} \mathcal{U} \times \mathcal{G}_c(\mathsf{L}_0) & \to & \pi^{-1}(\mathcal{U}) \\ ([\mathsf{L}], \mathsf{L}') & \mapsto & s([\mathsf{L}])(\mathsf{L}') \end{array}$$

is a local trivialization, proving the result. ∎

By theorem 2.3, the long exact sequence in homotopy and the contractibility of $\Lambda$, we have the next result.

Corollary 2.4 *There exists isomorphisms* $\pi_j(\Lambda_e) \approx \pi_{j-1}(\mathcal{G}_c(\mathsf{L}_0)) \approx \pi_{j-1}(\mathsf{U}(\infty)/\mathsf{O}(\infty))$ *for all* $j \in \mathbb{N}$.

A pair $(\mathsf{L}_0, \mathsf{L}_1)$ of lagrangians of $\mathsf{H}$ is said *complementary* if $\mathsf{L}_0 \cap \mathsf{L}_1 = 0$ and $\mathsf{L}_0 + \mathsf{L}_1 = \mathsf{H}$. The set of all lagrangians complementary to $\mathsf{L}_1$ is denoted by $\Lambda_0(\mathsf{L}_1)$. Define

$$\begin{array}{rcl} \varphi_{\mathsf{L}_0,\mathsf{L}_1} : \Lambda_0(\mathsf{L}_1) & \to & \mathscr{L}^{\mathrm{sa}}(\mathsf{L}_0) \\ \mathsf{L} & \mapsto & P_{\mathsf{L}_0} J S, \end{array}$$

where $S : \mathsf{L}_0 \to \mathsf{L}_1$ has graph $\mathrm{Gr}(S) = \mathsf{L}$ and $\mathscr{L}^{\mathrm{sa}}(\mathsf{L}_0)$ denotes the Banach space of bounded self-adjoint operators in $\mathsf{L}_0$. This map is a bijection and the collection $\{\varphi_{\mathsf{L}_0,\mathsf{L}_1}\}$, where $(\mathsf{L}_0, \mathsf{L}_1)$ runs over the set of all complementary lagrangians in $\mathsf{H}$, is a real-analytic atlas in $\Lambda$ which makes it a Banach manifold modeled on the space of bounded self-adjoint operators on a Hilbert space. More details can be found in [3], [6].

Lemma 2.5 *Given a pair* $(\mathsf{L}_0, \mathsf{L}_1)$ *of complementary lagrangians in* $\mathsf{H}$ *and* $\mathsf{L} \in \Lambda_0(\mathsf{L}_1)$*, we have that* $(\mathsf{L}, \mathsf{L}_0)$ *is a Fredholm pair if and only if* $\varphi_{\mathsf{L}_0,\mathsf{L}_1}(\mathsf{L})$ *is a Fredholm operator.*

Proof. Since $\ker \varphi_{\mathsf{L}_0,\mathsf{L}_1}(\mathsf{L}) = \mathsf{L} \cap \mathsf{L}_0$, we have that $\dim\ker \varphi_{\mathsf{L}_0,\mathsf{L}_1}(\mathsf{L}) = \dim(\mathsf{L} \cap \mathsf{L}_0)$. The inequalities

$$\gamma(\mathsf{L}_0, \mathsf{L}_1)\,\mathrm{dist}\,(u, \mathsf{L} \cap \mathsf{L}_0) \le \mathrm{dist}\,(u + Su, \mathsf{L} \cap \mathsf{L}_0) \le \|I + S\|\,\mathrm{dist}\,(u, \mathsf{L} \cap \mathsf{L}_0)\,,$$

where $\gamma(\mathsf{L}_0, \mathsf{L}_1)$ is the usual minimum gap between $\mathsf{L}_0$ and $\mathsf{L}_1$ (see [5, Chapter IV]) and $S : \mathsf{L}_0 \to \mathsf{L}_1$ has graph $\mathsf{L}$, imply that

$$\|I + S\|^{-1}\frac{\|P_{\mathsf{L}_0}JSu\|}{\mathrm{dist}\,(u, \mathsf{L} \cap \mathsf{L}_0)} \le \frac{\mathrm{dist}\,(u + Su, \mathsf{L}_0)}{\mathrm{dist}\,(u + Su, \mathsf{L} \cap \mathsf{L}_0)} \le \gamma(\mathsf{L}_0, \mathsf{L}_1)^{-1}\frac{\|P_{\mathsf{L}_0}JSu\|}{\mathrm{dist}\,(u, \mathsf{L} \cap \mathsf{L}_0)}\,,$$

so,

$$\|I + S\|^{-1}\gamma(P_{\mathsf{L}_0}JS) \le \gamma(\mathsf{L}, \mathsf{L}_0) \le \gamma(\mathsf{L}_0, \mathsf{L}_1)^{-1}\gamma(P_{\mathsf{L}_0}JS)\,,$$

where $\gamma(P_{\mathsf{L}_0}JS)$ denotes the minimum modulus of $P_{\mathsf{L}_0}JS$, as defined in [5, Chapter IV]. In particular, $\mathsf{L} + \mathsf{L}_0$ is closed if and only if $P_{\mathsf{L}_0}JS$ has closed range. Since $P_{\mathsf{L}_0}JS$ is self-adjoint and $(\mathsf{L} + \mathsf{L}_0)^{\perp} = J(\mathsf{L}_0 \cap \mathsf{L}_1)$, we conclude the result. ∎

Since the set of Fredholm self-adjoint operators is open in the space of all self-adjoint operators, we conclude that $\mathscr{F}_{\mathsf{L}_0}(\Lambda) \subset \Lambda$ is open. Before proving our principal result, we need a simple lemma which says that a Fredholm pair of lagrangians is a complementary pair modulo a compact perturbation.

Lemma 2.6 *If* $(\mathsf{L}_0, \mathsf{L}_1)$ *is a Fredholm pair of lagrangians then there exists* $\mathsf{L}_1' \in \Lambda$, $\mathsf{L}_1' \sim \mathsf{L}_1$*, such that* $(\mathsf{L}_0, \mathsf{L}_1')$ *is a complementary pair.*

Proof. Let $\mathsf{L}_2$ be a lagrangian complementary to $\mathsf{L}_1$ and $\mathsf{L}_0$ and $S \in \mathscr{L}(\mathsf{L}_0, \mathsf{L}_2)$ such that $\mathrm{Gr}(S) = \mathsf{L}_1$. By lemma 2.5, $T = \varphi_{\mathsf{L}_0,\mathsf{L}_2}(\mathsf{L}_1)$ is a Fredholm operator in $\mathsf{L}_0$. Defining $T' = T + P_{\ker T}$ and $S' = (P_{\mathsf{L}_0}J|_{\mathsf{L}_2})^{-1}T'$, we obtain that $T'$ is a isomorphism and $\mathsf{L}_1' = \mathrm{Gr}(S')$ is a lagrangian complementary to $\mathsf{L}_0$. Since $S'$ is a compact perturbation of $S$, we have that $\mathsf{L}_1' = \mathrm{Gr}(S') \sim \mathrm{Gr}(S) = \mathsf{L}_1$. This proves the result. ∎

The next theorem is the principal result of this article.

Theorem 2.7 *There exists a homotopy equivalence* $\mathscr{F}_{\mathsf{L}_0}(\Lambda) \approx \mathrm{U}(\infty)/\mathrm{O}(\infty)$.

Proof. Consider the restriction $\chi = \pi|_{\mathscr{F}_{\mathsf{L}_0}(\Lambda)} : \mathscr{F}_{\mathsf{L}_0}(\Lambda) \to \pi(\mathscr{F}_{\mathsf{L}_0}(\Lambda))$ of the fiber bundle $\pi : \Lambda \to \Lambda_e$. Since $\mathscr{F}_{\mathsf{L}_0}(\Lambda)$ is stable under compact perturbations, the map $\chi$ defines a fiber bundle with typical fiber $\mathcal{G}_c(\mathsf{L}_0)$. The previous lemma shows that $\pi(\mathscr{F}_{\mathsf{L}_0}(\Lambda)) = \pi(\Lambda_0(\mathsf{L}_0))$. The set $\Lambda_0(\mathsf{L}_0)$ is homeomorphic to $\mathscr{L}^{\mathrm{sa}}(\mathsf{L}_0^{\perp})$ via the usual map

$$\begin{aligned} \varphi_{\mathsf{L}_0^{\perp},\mathsf{L}_0} : \Lambda_0(\mathsf{L}_0) &\to \mathscr{L}^{\mathrm{sa}}(\mathsf{L}_0^{\perp}) \\ \mathsf{L} &\mapsto JS\,, \end{aligned}$$

where $S : \mathsf{L}_0^{\perp} \to \mathsf{L}_0$ has graph $\mathrm{Gr}(S) = \mathsf{L}$ and $\mathscr{L}^{\mathrm{sa}}(\mathsf{L}_0^{\perp})$ denotes the Banach space of bounded self-adjoint operators in $\mathsf{L}_0$. Now, the map

$$\begin{aligned} \pi(\Lambda_0(\mathsf{L}_0)) \times [0,1] &\to \pi(\Lambda_0(\mathsf{L}_0)) \\ (\pi(\varphi^{-1}_{\mathsf{L}_0^{\perp},\mathsf{L}_0}(T)), t) &\mapsto \pi(\varphi^{-1}_{\mathsf{L}_0^{\perp},\mathsf{L}_0}(tT)) \end{aligned}$$

is well-defined, continuous and gives a homotopy between the identity function on $\pi(\Lambda_0(\mathsf{L}_0))$ and a constant. So, $\pi(\Lambda_0(\mathsf{L}_0))$ is contractible, and the long exact sequence of the fibration applied to $\chi : \mathscr{F}_{\mathsf{L}_0}(\Lambda) \to \pi(\Lambda_0(\mathsf{L}_0))$ implies that exists a weak homotopy equivalence between $\mathscr{F}_{\mathsf{L}_0}(\Lambda)$ and $\mathcal{G}_c(\mathsf{L}_0)$. The theorem follows by [8, Theorem 15] and proposition 2.2. ∎

We end with a result about the Maslov index for curves in $\mathscr{F}_{\mathsf{L}_0}(\Lambda)$. Given a curve $\gamma : [0,1] \to \Lambda_0(\mathsf{L}_1)$, the Maslov index of $\gamma$ with respect to $\mathsf{L}_0$ is defined as

$$\mu_{\mathsf{L}_0}(\gamma) = \mathrm{sf}(\varphi_{\mathsf{L}_0,\mathsf{L}_1} \circ \gamma)\,,$$

where $\mathrm{sf}(\varphi_{\mathsf{L}_0,\mathsf{L}_1} \circ \gamma)$ denotes the spectral flow of the curve of Fredholm self-adjoint operators $\varphi_{\mathsf{L}_0,\mathsf{L}_1} \circ \gamma$, as defined in [9]. A staightforward calculation show that this definition is coherent with the open covering $\{\Lambda_0(\mathsf{L}_1)\}_{\mathsf{L}_1 \in \Lambda}$, that is, if $\gamma$ has image in $\Lambda_0(\mathsf{L}_1) \cap \Lambda_0(\mathsf{L}_2)$, then $\mathrm{sf}(\varphi_{\mathsf{L}_0,\mathsf{L}_1} \circ \gamma) = \mathrm{sf}(\varphi_{\mathsf{L}_0,\mathsf{L}_2} \circ \gamma)$. The integer $\mu_{\mathsf{L}_0}(\gamma)$ is additive by concatenation and invariant by a homotopy fixing the endpoints. Given any curve $\gamma : [0,1] \to \mathscr{F}_{\mathsf{L}_0}(\Lambda)$, we obtain a partition $0 = x_0 < \ldots < x_n = 1$ and $\mathsf{L}_1, \ldots, \mathsf{L}_n \in \Lambda$ such that $(\mathsf{L}_0, \mathsf{L}_j)$ are complementary and $\gamma_j = \gamma|_{[x_{j-1},x_j]}$ has image in $\Lambda_0(\mathsf{L}_j)$, for $j = 1, \ldots, n$. We define the Maslov index of $\gamma$ with respect to $\mathsf{L}_0$ as

$$\mu_{\mathsf{L}_0}(\gamma) = \sum_{j=1}^{n} \mathrm{sf}(\varphi_{\mathsf{L}_0,\mathsf{L}_j} \circ \gamma_j)\,.$$

Homotopy invariance and additivity properties show that the Maslov index is well-defined and gives a homomorphism of the fundamental grupoid[1] of $\mathscr{F}_{\mathsf{L}_0}(\Lambda)$ with $\mathbb{Z}$. For more details about the construction of the Maslov index, see [3], [4], [6]. The next proposition shows that the restriction of $\mu_{\mathsf{L}_0}$ to the fundamental group of $\mathscr{F}_{\mathsf{L}_0}(\Lambda)$ is a isomorphism.

Proposition 2.8 *The Maslov index with respect to a fixed lagrangian* $\mathsf{L}_0$ *is a isomorphism* $\mu_{\mathsf{L}_0} : \pi_1(\mathscr{F}_{\mathsf{L}_0}(\Lambda)) \to \mathbb{Z}$.

Proof. Note that by theorem 2.7, $\mathscr{F}_{\mathsf{L}_0}(\Lambda)$ is path-connected and $\pi_1(\mathscr{F}_{\mathsf{L}_0}(\Lambda)) = \mathbb{Z}$, so it is sufficient to prove that $\mu_{\mathsf{L}_0}$ is onto. Given $n \in \mathbb{Z}$, take a loop $\alpha$ in $\mathscr{F}^{\mathrm{sa}}(\mathsf{L}_0^{\perp})$ with spectral flow $n$; the curve $\gamma = \varphi^{-1}_{\mathsf{L}_0^{\perp},\mathsf{L}_0} \circ \alpha$ is a loop in $\mathscr{F}_{\mathsf{L}_0}(\Lambda)$ with Maslov index $n$, as desired. ∎

[1] The fundamental grupoid of a space $X$ is the set of all homotopy classes, with fixed endpoints, of curves in $X$.

Departamento de Matemática,
Instituto de Matemática e Estatística
Universidade de São Paulo, Brasil.
E-mail addresses: zieca@ime.usp.br, piccione@ime.usp.br